# White-noise-aided Control

Wuhua Hu

School of EEE, Nanyang Technological University, Singapore

E-mail: haitun_ahua@yahoo.com.cn

***Abstract:*** The issue of white-noise-aided control is considered and its availability is proved. And a noise-aiding way is developed to stabilize perturbed systems to be input-to-state stable (ISS) with respect to (w.r.t.) perturbations. To illustrate its effectiveness, the white-noise-aided control of a parameter perturbed chaotic Chen system is given as an example. And numerically, it shows that, comparing to the un-noise-aided case, noise-aided control can not only shorten the control's transient process but also save its cost. These are also demonstrated by various aiding noises such as common (symmetric) noise and non-common (independent or asymmetric) noise, where common noise is found to be the most efficient in enhancing the control.

Keywords: White noise; Control; Input-to-state stability; Robust stability; Perturbation; Chaos

## 1. Introduction

Noise has been found to act positively in various areas such as physics, chemistry, biology and ecology [1-5]. However, functions of the noise-enhanced and the noise-induced effects have mainly been studied numerically and experimentally [6-9], and little theoretical analysis has been done yet. For the latter, the creative work of [10] is commendable, which uncover the reason for noise-induced synchronization. And this could also be used to explain many observed noise-enhanced effects, such as synchronization [7;9] and stability [3;11]. Besides, we should notice the pioneer work of [8], which numerically studied the constructive role of external noise in the control of chaotic dynamics. These all urge us to pay more attention to the theoretical explanation of the positive effects of noise.

For anther thing, it was discovered that when realizing the same synchronization, the intensities of needed noise is much weaker than the coupling strengths of an alternative linear feedback [10]. And this is also supported by the facts that even much small noise can enhance stability greatly [3], or push the system to synchronize more [9]. Moreover, the paper [8] finds that noise can function to shorten the transient process in a control task. These benefits all evoke us to make use of noise, e.g., to save the energy in doing a control task, rather than simply suppress them in a traditional sense (e.g., [12]).

The above two points contribute to our idea that making use of the constructive role of noise to improve a job. And here we focus this job on control, which is actually a very general issue. More specifically, we are to study how the noise can aid to stabilize a system to be ISS w.r.t. the external input (also called "controls" or "disturbances" depending on the context). The ISS case is considered here, because it has been recently developed to cover much general and common stability or control tasks [13;14].

In the following context, we will mainly prove the availability of noise-aided control, with a sufficient condition for such control given, and we will illustrate the way of designing such control with an example of stabilizing a chaotic Chen system [15] to be ISS w.r.t. external disturbance (including white and non-white noise ones). The energy cost of such control and its dynamic response are also compared with those un-noise-aided cases. Also, the impacts of noises with different coherences on the control effect are considered.

## 2. Preliminaries on ISS

**Notations:** 1) A $\mathcal{K}$-function means a function $\alpha : [0,a) \to [0,\infty)$ that is continuous and strictly increasing and satisfies $\alpha(0) = 0$; furthermore, if $a = \infty$ and $\alpha(r) \to \infty$ when $r \to \infty$, then $\alpha$ is called a $\mathcal{K}_\infty$-function; 2) A $\mathcal{KL}$-function means a function $\beta : [0,a) \times [0,\infty) \to [0,\infty)$ that is a class of $\mathcal{K}$ on the first argument and satisfies $\beta(r,s) \to 0$ as $s \to \infty$; 3) $\|\bullet\|$ denotes the Euclidean norm, and $\|\bullet\|_\infty$ denotes the supremum norm; 4) $\mathcal{M}_{\mathcal{O}} :=$ the set of all measurable functions from $\mathbb{R}_{\geq 0}$ to the unit ball in $\mathbb{R}^p$. 5) Matrices "$> 0$", "$\geq 0$", or "$< 0$" respectively means it is positive definite, positive semi-definite and negative definite.

Consider the continuous time system of the standard form

$$\dot{x} = f(x, u), \tag{1}$$

where $x(t) \in \mathbb{R}^n$ and $u(t) \in \mathbb{R}^m$. It is assumed that $f : \mathbb{R}^n \times \mathbb{R}^m$ is locally Lipschitz and satisfies $f(0,0) = 0$. Controls or inputs are measurable locally essentially bounded functions $u : \mathbb{R}_{\geq 0} \times \mathbb{R}^m$. The set of all such functions is denoted by $L_\infty^m$, $x(t, x_0, u)$ denotes the trajectory of system (1) with initial state $x(0) = x_0$ and input $u$. This is a priori defined only on some maximal interval $[0, T_{x_0,u})$, with $T_{x_0,u} \leq +\infty$. The system is (forward-) complete if $T_{x_0,u} = +\infty$ for all $x_0$ and $u$.

**Definition 1.** [14;16] System (1) is said to be *input-to-state stable (ISS)* if there exists a $\mathcal{KL}$-function $\beta$ and a $\mathcal{K}$-function $\gamma$ such that, for any input $u \in L_\infty^m$ and any $x_0 \in \mathbb{R}^n$, it holds that

$$\|x(t, x_0, u)\| \leq \beta(\|x_0\|, t) + \gamma(\|u\|_\infty), \ \forall t \geq 0. \tag{2}$$ □

**Definition 2.** [14;16] System (1) is said to be *robustly stable (RS)* if there exists a $\mathcal{K}_\infty$-function $\rho$ such that the system

$$x = g(x, d) := f(x, d\rho(\|x\|)) \tag{3}$$

is globally asymptotically stable (GAS) uniformly in this sense: for some $\mathcal{KL}$-function $\beta$, the estimate

$$\|x(t, x_0, u)\| \leq \beta(\|x_0\|, t), \ \forall t \geq 0, \tag{4}$$

holds for all $d \in \mathcal{M}_{\mathcal{O}}$ and any $x_0 \in \mathbb{R}^n$. □

**Lemma 1.** [14;16] *System (1) is ISS if and only if (iff) it is RS.* □

The following definition and lemma are derived by reducing those for ISS to non-disturbed cases (see [17]).

**Definition 3.** A smooth function $V : \mathbb{R}^n \to \mathbb{R}_{\geq 0}$ is called an *GAS-control Lyapunov function (GAS-CLF)* for system (1), if there exist $\mathcal{K}_\infty$-functions $\alpha_1, \alpha_2$, such that

$$\alpha_1(\| x \|) \leq V(x) \leq \alpha_2(\| x \|), \text{ and } \nabla V(x) f(x,u) < 0$$

for all $x \neq 0$ and $u \in \mathbb{R}^m$. □

**Lemma 2.** *System (1) is GAS stabilizable iff there exists a GAS-CLF.* □

If define $u := d\rho(\| x \|)$, Lemma 1 indicates that system (3) is uniformly GAS (UGAS) stabilizable iff there exists a GAS-CLF.

Next, we consider the system affine in control

$$x = f(x) + G(x)[k(x) + u], \ f(0) = 0, \tag{5}$$

and the following definition is made:

**Definition 4.** [18] System (5) is *smoothly stabilizable* if there exists a smooth map $k : \mathbb{R}^n \to \mathbb{R}^m$ with $k(0) = 0$ such that (5) with $u \equiv 0$ is GAS. It is *smoothly ISS stabilizable* if there is such a $K$ so that system (5) becomes ISS w.r.t. $u$.

And a conclusion can be obtained.

**Lemma 3.** [14;18] *Smooth stabilizability implies smooth input to state stabilizability.*

Note that somewhat less than smooth (differentiability) of $k$ is enough for the above argument: continuity is enough. However, if no continuous feedback stabilizer exists, then no smooth CLF $V$ can be found (Continuous stabilization of nonlinear systems is basically equivalent to the existence of what are called smooth CLF.) [14].

## 3. Main results

Let the nominal system be

$$\dot{x} = A_0 x + f_0(x), \tag{6}$$

where constant matrix $A_0 \in \mathbb{R}^{n\times n}$ and the nonlinear vector function $f_0 : \mathbb{R}^n \to \mathbb{R}^n$, which is smooth with $f(0) = 0$. Let $u : \mathbb{R}_{\geq 0} \to \mathbb{R}^m$ be the controls, and $\dot{B}_c : \mathbb{R}_{\geq 0} \to \mathbb{R}^l$ be the white noises to aid the control, and $\omega : \mathbb{R}_{\geq 0} \to \mathbb{R}^p$ , $\dot{B}_d : \mathbb{R}_{\geq 0} \to \mathbb{R}^p$ be respectively the non-white and white noise disturbances (e.g., the combined model uncertainties and actuator disturbances). Besides, suppose $\dot{B}_{c,d}$ have mean-values of $0$ and covariances of $1$, where $B_{c,d}$ are the standard Wiener processes [19]. We consider the control system on $\mathbb{R}^n$ of the general form (see Fig. 1):

$$\begin{aligned}\dot{x} &= A_0 x + g(x, u, \sigma_c \dot{B}_c, \sigma_d \dot{B}_d + \omega) \\ &:= A_0 x + f_0(x) + G(x)u + C(x)\sigma_c \dot{B}_c + D(x)(\sigma_d \dot{B}_d + \omega),\end{aligned} \tag{7}$$

where $\sigma_c = diag(\sigma_1^c, \sigma_2^c, \cdots, \sigma_l^c)$ , $\sigma_d = diag(\sigma_1^d, \sigma_2^d, \cdots, \sigma_p^d)$ are respectively the intensity matrices of $\dot{B}_c$, $\dot{B}_d$, and $C = (C_1, C_2, \cdots, C_l) : \mathbb{R}^n \to \mathbb{R}^{n\times l}$, $D = (D_1, D_2, \cdots, D_p) : \mathbb{R}^n \to \mathbb{R}^{n\times p}$ and $G : \mathbb{R}^n \to \mathbb{R}^{n\times m}$.

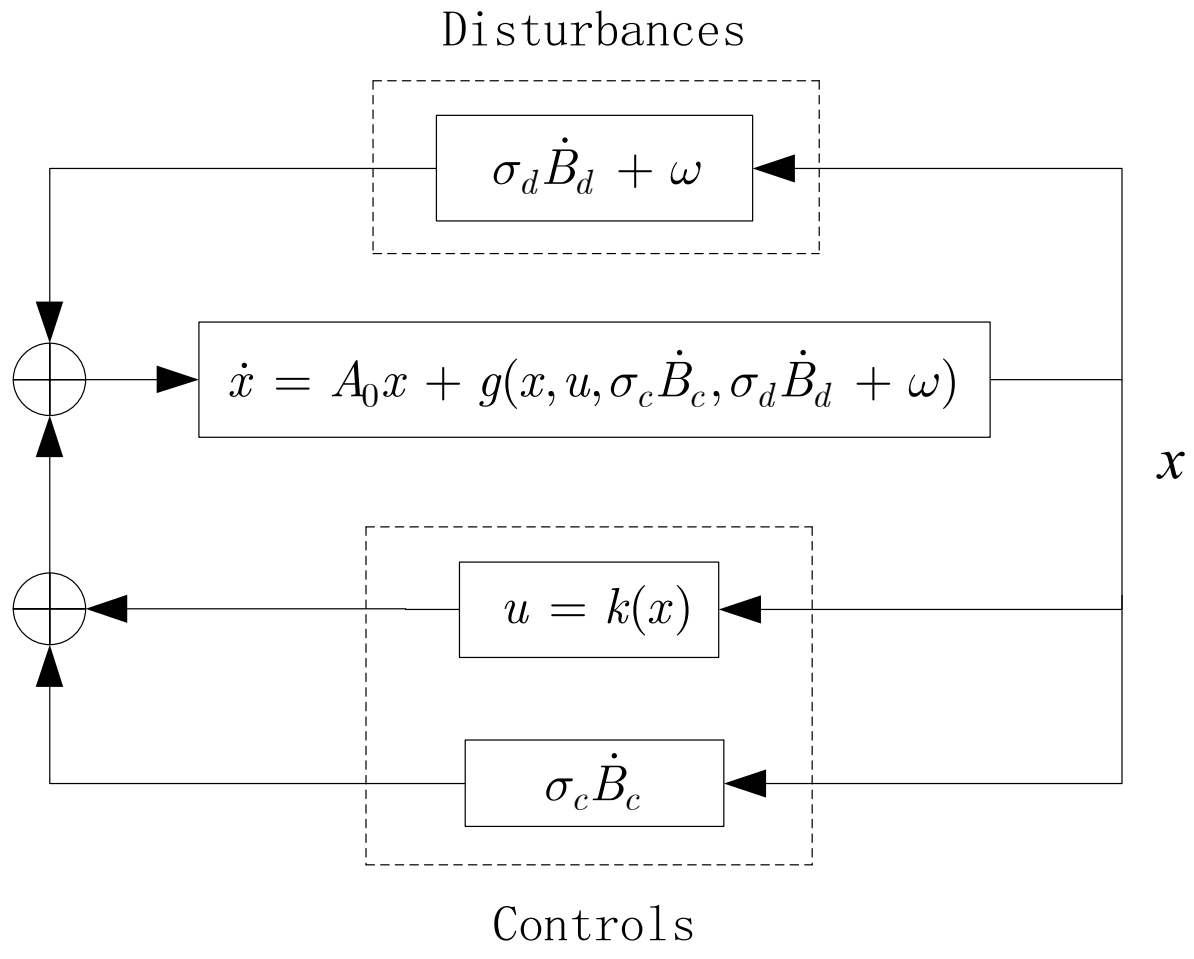

Fig. 1 Closed-loop system (7)

(**Declaration:** *in this paper, since the disturbance containing white noise may be dealt with, the*

*stability referred and deduction done, in the presentation of noise, all mean in an almost sure sense; and for simplicity, we omit the words "almost surely" throughout the paper.*)

Then, we are to consider how this system can be stabilized ISS w.r.t. the disturbance $\omega$, under some feedback law $k(x)$ and some aiding noise $\dot{B}_c$. First, rewrite (7) as

$$\dot{x} = A_0 x + f_0(x) + G(x)k(x) + H(x)\dot{B} + D(x)\omega, \tag{8}$$

where $H(x) := [C(x)\sigma_c \;\; D(x)\sigma_d] \in \mathbb{R}^{n\times(l+p)}$ and $\dot{B} := \left[\dot{B}_c^T \;\; \dot{B}_d^T\right]^T \in \mathbb{R}^{l+p}$. Suppose $G(x)k(x) = \tilde{A}_0 x + \tilde{f}_0(x)$. Then let $A := A_0 + \tilde{A}_0$ and $f(x) := f_0(x) + \tilde{f}_0(x)$, and so

$$\dot{x} = Ax + f(x) + H(x)\dot{B} + D(x)\omega. \tag{9}$$

Hence, our task is now to find some proper white noises $\dot{B}$ (i.e., $\dot{B}_c$) to guarantee (9) to be ISS w.r.t. $\omega$. And we make **several assumptions:** A1) *The solution of system (7) uniquely exists;* A2) *The components of* $\dot{B}$ *are independent of each other;* A3) *System (7) is smoothly stabilizable when* $\omega \equiv 0$*;* A4) $f(x)$ *satisfies the Lipschitz condition*

$$f^T(x)f(x) \le x^T L x. \tag{10}$$

A3 implies that (7) is smoothly ISS stabilizable w.r.t. $\omega$ (Lemma 3). Further, due to the equivalence of ISS and RS (Lemma 1), if (9) were ISS, then there should exist a $\mathcal{K}_\infty$-function $\rho$ such that the system

$$\dot{x} = Ax + f(x) + H(x)\dot{B} + \rho(\|x\|)D(x)d \tag{11}$$

is UGAS w.r.t. all $d \in \mathcal{M}_{\mathcal{O}}$. Then, similar to the theorem in [10], a sufficient condition for such a goal can be obtained.

For the convenience of statement, we first give an auxiliary lemma.

**Lemma 4.** Given $R(\in \mathbb{R}^{n\times n}) \ge 0$ and $Y : \mathbb{R}^n \to \mathbb{R}^n$, for $\forall x \in \mathbb{R}^n$ if $\exists S(\in \mathbb{R}^{n\times n}) \ge 0$ s.t. $Y(x)^T R Y(x) \le x^T S x$, then $\exists T(\in \mathbb{R}^{n\times n}) \ge 0$ s.t. $x^T R Y(x) \le x^T T x$. □

*Proof.* Since $R \ge 0$, there exists a symmetric matrix $B$ satisfying $R = B^2$. So,

$$x^T R Y(x) = x^T B^2 Y(x) \le \left\| x^T B \right\| \left\| BY(x) \right\| = \sqrt{x^T BBx}\sqrt{Y^T(x)BBY(x)}$$
$$= \sqrt{x^T Rx}\sqrt{Y^T(x)RY(x)} \le \sqrt{x^T Rx}\sqrt{x^T Sx} \le \sqrt{rs}\, \| x \|^2 ,$$

where $r \ge 0,\ s \ge 0$ are the largest eigenvalues of $R,\ S$ respectively. Let $T := \sqrt{rs} I_n \in \mathbb{R}^{n\times n}$, with $I_n$ being the unite matrix in $\mathbb{R}^n$. Hence, $T$ is positive semi-definite and $x^T R Y(x) \le x^T T x$ holds. □

**Theorem 1.** *System (11) (and thus (9) and (7)) is ISS w.r.t.* $\omega$*, if, besides the assumptions A1, A2, A3 and A4, there exists a symmetric matrix* $P > 0$ *and* $R \ge 0,\ K_i^c \ge 0,\ J_i^c \ge 0,\ K_j^d \ge 0,\ J_j^d \ge 0\ \ (i = 1,2,\cdots,l;\ j = 1,2,\cdots,p)$*, satisfying all the three conditions*

$$a)\quad \begin{array}{c} \rho(\| x \|)\left\| x^T PD(x) \right\| \le x^T Rx, \\ x^T PD_j(x) \ge x^T J_j^d x,\ D_j^T(x)PD_j(x) \le x^T K_j^d x, \end{array} \qquad (12)$$

$$b)\quad x^T PC_i(x) \ge x^T J_i^c x,\ C_i^T(x)PC_i(x) \le x^T K_i^c x, \qquad (13)$$

$$c)\quad \begin{array}{l} A^T P + PA + \varepsilon P^2 + \dfrac{1}{\varepsilon} L + 2R + \sum_{i=1}^{l} \sigma_i^c K_i^c + \sum_{i=1}^{p} \sigma_i^d K_i^d \\ -\dfrac{2}{\beta_{\max}}\left( \sum_{i=1}^{l} (\sigma_i^c)^2 \alpha_i^c J_i^c + \sum_{i=1}^{p} (\sigma_i^c)^2 \alpha_i^d J_i^d \right) < 0, \end{array} \qquad (14)$$

*where* $\alpha_i^c,\ \alpha_i^d$ *are respectively the smallest eigenvalues of* $C_i,\ D_i$*, and* $\beta_{\max}$ *is the largest eigenvalue of* $P$*.* □

*Proof.* We precede the proof similar to that in [10]. Due to the nonzero property of $x(t;x_0)$ when $x_0 \ne 0$ (see Appendix A in [10]), we can introduce the following function:

$$V(x) = \frac{1}{2}\log x^T Px = \frac{1}{2}\log\left\| P^{1/2} x \right\|. \qquad (15)$$

By applying Itô's formula to (15) along with system (11), it yields that

$$V(x(t)) = V(x(t_0)) + \int_{t_0}^{t} \left\{ \begin{array}{l} \dfrac{\partial V(x(s))}{\partial x(s)}[Ax(s) + f(x(s)) + \varphi(x(s))D(x(s))d] \\ +\dfrac{1}{2}\mathrm{trace} H^T(x(s))\dfrac{\partial^2 V(x(s))}{\partial x^T(s)\partial x(s)} H(x(s)) \end{array} \right\} ds + M(t), (16)$$

where

$$\frac{\partial V(x(s))}{\partial x(s)} = \frac{x(s)^T P}{\left\| P^{1/2} x(s) \right\|^2},$$

$$\frac{\partial^2 V(x(s))}{\partial x^T(s) \partial x(s)} = \frac{P}{\left\| P^{1/2} x(s) \right\|^2} - 2\frac{P x(s) x^T(s) P}{\left\| P^{1/2} x(s) \right\|^4}, \tag{17}$$

and the continuous martingale is

$$M(t) = \int_{t_0}^{t} \frac{\partial V(x(s))}{\partial x(s)} H(x(s)) dB(s) = \int_{t_0}^{t} \frac{x(s)^T P H(x(s))}{\left\| P^{1/2} x(s) \right\|^2} dB(s), \tag{18}$$

with $M(t_0) = 0$, and the quadratic variation is

$$[M(t), M(t)] = \int_{t_0}^{t} \frac{\left\| x(s)^T P H(x(s)) \right\|^2}{\left\| P^{1/2} x(s) \right\|^4} ds. \tag{19}$$

Hence,

$$\begin{aligned}
[M(t), M(t)] &= \int_{t_0}^{t} \frac{\left\| x(s)^T P C(x(s)) \sigma_c \right\|^2 + \left\| x(s)^T P D(x(s)) \sigma_d \right\|^2}{\left\| P^{1/2} x(s) \right\|^4} ds \\
&= \int_{t_0}^{t} \frac{\sum_{i=1}^{l} \left[ \sigma_i^c x(s)^T P C_i(x(s)) \right]^2 + \sum_{i=1}^{p} \left[ \sigma_i^d x(s)^T P D_i(x(s)) \right]^2}{\left\| P^{1/2} x(s) \right\|^4} ds \\
&\geq \int_{t_0}^{t} \frac{\sum_{i=1}^{l} \left[ \sigma_i^c x(s)^T J_i^c x(s) \right]^2 + \sum_{i=1}^{p} \left[ \sigma_i^d x(s)^T J_i^d x(s) \right]^2}{\left\| P^{1/2} x(s) \right\|^4} ds \\
&\geq c(t - t_0),
\end{aligned}$$

with

$$c = \frac{\sum_{i=1}^{l} (\sigma_i^c \eta_i^c)^2 + \sum_{i=1}^{p} (\sigma_i^d \eta_i^d)^2}{\beta_{\max}^2},$$

where $\eta_i^c, \eta_i^d$ are the smallest eigenvalues of $J_i^c$, $J_i^d$ respectively, and $\beta_{\max}$ is the largest eigenvalue of $P$. And according to Lemma 4, it follows that

$$[M(t), M(t)] \leq \int_{t_0}^{t} \frac{\sum_{i=1}^{l} \left[ \sigma_i^c x(s)^T J_i'^c x(s) \right]^2 + \sum_{i=1}^{p} \left[ \sigma_i^d x(s)^T J_i'^d x(s) \right]^2}{\left\| P^{1/2} x(s) \right\|^4} ds \leq c'(t - t_0),$$

with

$$c' = \frac{\sum_{i=1}^{l} \left(\sigma_i^c \eta_i'^c\right)^2 + \sum_{i=1}^{p} \left(\sigma_i^d \eta_i'^d\right)^2}{\beta_{\min}^2},$$

where $\eta_i'^c$, $\eta_i'^d$ are the largest eigenvalues of $J_i'^c$, $J_i'^d$ respectively, and $\beta_{\min}$ is the smallest eigenvalue of $P$.

So, $c(t-t_0) \le [M(t), M(t)] \le c'(t-t_0)$ for some constants $c, c' > 0$ for all $t \ge t_0$. And this implies (Theorem 3.29 and Theorem 7.33 in [19]) that

$$\lim_{x\to\infty} \frac{M(t)}{t} = 0. \tag{20}$$

Moreover, from (10) and an elementary inequality for vectors, it follows that for any $\varepsilon > 0$,

$$x^T P f(x) \le x^T \left(\frac{\varepsilon}{2} P^2 + \frac{1}{2\varepsilon} L\right) x. \tag{21}$$

Now, substituting inequalities (13) and (21) into (16) gives

$$\begin{aligned}
\frac{V(x(t))}{t} &\le \frac{V(x(0))}{t} + \frac{1}{t}\int_{t_0}^{t} x^T(s)\left(\frac{A^T P + PA}{2} + \frac{\varepsilon}{2}P^2 + \frac{1}{2\varepsilon}L + R\right)\frac{x(s)}{\left\|P^{1/2}x(s)\right\|^2}\,ds \\
&\quad + \frac{1}{t}\int_{t_0}^{t} \operatorname{trace}\left[\frac{H^T P H}{2\left\|P^{1/2}x(s)\right\|^2} - \frac{H^T P x(s) x^T(s) P H}{\left\|P^{1/2}x(s)\right\|^4}\right] ds + \frac{M(s)}{t} \\
&= \frac{V(x(0))}{t} + \frac{M(s)}{t} + \frac{1}{t}\int_{t_0}^{t} x^T(s)\left(\frac{A^T P + PA}{2} + \frac{\varepsilon}{2}P^2 + \frac{1}{2\varepsilon}L + R\right)\frac{x(s)}{\left\|P^{1/2}x(s)\right\|^2}\,ds \\
&\quad + \frac{1}{t}\int_{t_0}^{t} \left\{\begin{aligned} &\sum_{i=1}^{l}\left\{\frac{\sigma_i^c C_i^T P C_i}{2\left\|P^{1/2}x(s)\right\|^2} - \frac{\left[\sigma_i^c x^T(s) P C_i\right]^2}{\left\|P^{1/2}x(s)\right\|^4}\right\} \\ &+ \sum_{i=1}^{p}\left\{\frac{\sigma_i^d D_i^T P D_i}{2\left\|P^{1/2}x(s)\right\|^2} - \frac{\left[\sigma_i^d x^T(s) P D_i\right]^2}{\left\|P^{1/2}x(s)\right\|^4}\right\}\end{aligned}\right\} ds \\
&\le \frac{V(x(0))}{t} + \frac{M(s)}{t} - \frac{1}{2t}\int_{t_0}^{t} \frac{x^T(s) Q x(s)}{\left\|P^{1/2}x(s)\right\|^2} ds,
\end{aligned}$$

where

$$\begin{aligned}
-Q &= A^T P + PA + \varepsilon P^2 + \frac{1}{\varepsilon}L + 2R + \sum_{i=1}^{l}\sigma_i^c K_i^c + \sum_{i=1}^{p}\sigma_i^d K_i^d \\
&\quad - \frac{2}{\beta_{\max}}\left(\sum_{i=1}^{l}(\sigma_i^c)^2 \alpha_i^c J_i^c + \sum_{i=1}^{p}(\sigma_i^c)^2 \alpha_i^d J_i^d\right).
\end{aligned} \tag{22}$$

Since (14) holds, $Q > 0$. And because the assumption A3 implies the forward complete of system (7) (P. 9, [14]), the following inequalities would hold

$$\overline{\lim_{t\to\infty}} \frac{\log\left\|P^{1/2}x(t)\right\|}{t} = \overline{\lim_{t\to\infty}} \frac{V(x(t))}{t} \le \overline{\lim_{t\to\infty}} \left[ -\frac{1}{2t}\int_{t_0}^{t} \frac{x^T(s)Qx(s)}{\left\|P^{1/2}x(s)\right\|^2} ds \right] \le -\frac{\gamma}{2\beta_{\max}}, \quad (23)$$

with $\gamma > 0$ being the smallest eigenvalue of $Q$. (Here "lim sup" is required since the limit of $V(x(t))/t$ as $t \to \infty$ may not well exist). Therefore, (11) is UGAS with a robust margin of stability $\rho$. And hence (9) is RS and ISS w.r.t. $\omega$ (Lemma 1). So is (7). □

**Remark 1.** a) For linear systems, (21) changes to be $x^T Pf(x) = 0$, so the corresponding condition (14) can be modified by removing the term $\varepsilon P^2 + L/\varepsilon$. b) If $\omega \equiv 0$, then $R$ can be taken zero in (12), and system (7) would be GAS under the similar conditions. And if simultaneously $\dot{B}_d \equiv 0$, then the conditions for system (7) to be GAS simplify to be (13) together with

$$A^T P + PA + \varepsilon P^2 + \frac{1}{\varepsilon} L + \sum_{i=1}^{l} \sigma_i^c K_i^c - \frac{2}{\beta_{\max}} \sum_{i=1}^{l} (\sigma_i^c)^2 \alpha_i^c J_i^c < 0. \quad (24)$$

In this case, if we suppose $f_0 \equiv f$ namely $G(x)k(x) \equiv 0$, then it reduces to the situation that totally using noise to realize the control, and (13) together with (24) actually weaken the conditions obtained in [10]. □

Below, we give a corollary that asserts the probability one that white noise can be used to aid the control under some situation.

**Corollary 1.** *Let $C_i(x) = c_i x$ with $c_i > 0$, $i = 1, 2, \cdots, l$. Besides the assumptions A1, A2, A3 and A4, suppose there exist a symmetric matrix $P > 0$, with its eigenvalues satisfying $2\beta_{\min} > \beta_{\max}$, and such $R \ge 0$, $K_j^d \ge 0$, $J_j^d \ge 0$ $(j = 1, 2, \cdots, p)$ that (12) stands. Then, there would always exist such aiding noises $\sigma_c \dot{B}_c$ that system (11) is UGAS, and hence systems (9) and (7) are ISS w.r.t. $\omega$.* □

*Proof.* Since $C_i(x) = c_i x$, with $c_i > 0$ $(i = 1,2,\cdots,l)$, and there exists a symmetric $P > 0$ and such $R \geq 0$, $K_j^d \geq 0$ and $J_j^d \geq 0$ $(j = 1,2,\cdots,p)$ that (12) stand, we can let $J_i^c = c_i P > 0$, $K_i^c = c_i^2 P > 0$. And hence (13) holds. As a result, (22) can be obtained as

$$-Q := \left(A^T P + PA + \varepsilon P^2 + \frac{1}{\varepsilon} L + 2R + \sum_{i=1}^{p} \sigma_i^d K_i^d\right) - \frac{2}{\beta_{\max}} \sum_{i=1}^{p} (\sigma_i^d)^2 \alpha_i^d J_i^d - \left(\frac{2\beta_{\min}}{\beta_{\max}} - 1\right) \sum_{i=1}^{l} (\sigma_i^c c_i)^2 P.$$

Thus, as long as $\sigma_i^c$ is large enough $Q$ would be positive definite. And according to Theorem 1, the corresponding noises $\sigma_c \dot{B}_c$, together with the feedback law $k(x)$, would guarantee (11) to be UGAS, and hence (7) to be ISS w.r.t. $\omega$. □

With the case being considered that perturbations of white and non-white noises are not distinguishable or that there is no need to distinguish them, our task should change to design some feedback law $k(x)$ and some aiding noises $\sigma_c \dot{B}_c$ to stabilize the system

$$\dot{x} = Ax + f_0(x) + G(x)k(x) + C(x)\sigma_c \dot{B}_c + D(x)\tilde{\omega}, \tag{25}$$

to be ISS w.r.t. $\tilde{\omega}$, with $\tilde{\omega} := w + \sigma_d(t)\dot{B}_d$. By replacing $\dot{B}$ with $\dot{B}_c$ and $\omega$ with $\tilde{\omega}$, similar to Theorem 1, it can then be obtained that:

**Theorem 2.** *System (25) is ISS w.r.t.* $\tilde{\omega}$*, if, besides the assumptions A1, A2, A3 and A4, there exists a symmetric matrix* $P > 0$ *and* $R \geq 0$, $K_i^c \geq 0$, $J_i^c \geq 0$ $(i = 1,2,\cdots,l)$*, satisfying all the three conditions*

$$a)\ \rho(\|x\|)\left\|x^T PD(x)\right\| \leq x^T Rx, \tag{26}$$

$$b)\ x^T PC_i(x) \geq x^T J_i^c x,\ \ C_i^T(x)PC_i(x) \leq x^T K_i^c x, \tag{27}$$

$$c)\ A^T P + PA + \varepsilon P^2 + \frac{1}{\varepsilon} L + 2R + \sum_{i=1}^{l} \sigma_i^c K_i^c - \frac{2}{\beta_{\max}} \sum_{i=1}^{l} (\sigma_i^c)^2\, \alpha_i^c J_i^c < 0, \tag{28}$$

*where* $\alpha_i^c$ *is the smallest eigenvalue of* $C_i$*, and* $\beta_{\max}$ *is the largest eigenvalue of* $P$. □

*Proof.* The proof is almost the same as that of Theorem 1, so it is omitted. □

And the following corollary can be derived from Theorem 1.

**Corollary 2.** *Let* $C_i(x) = c_i x$ *with* $c_i > 0$*,* $i = 1, 2, \cdots, l$*. Besides the assumptions A1, A2, A3 and A4, suppose there exists a symmetric matrix* $P > 0$*, with its eigenvalues satisfying* $2\beta_{\min} > \beta_{\max}$*, and such* $R \geq 0$ *that* (26) *stands. Then, there would always exist such aiding noises* $\sigma_c \dot{B}_c$ *that system (25) is ISS w.r.t.* $\tilde{\omega}$. □

*Proof.* The proof is almost the same as that of Corollary 1, so it is omitted. □.

### 4. An example: stabilizing a perturbed chaotic Chen system to be ISS w.r.t. disturbances

For simplicity, suppose the disturbances are not distinguished and taken as a whole. Let us consider the parameters perturbed chaotic Chen system [15]

$$\begin{cases} \dot{x}_1 = (a - w_1)(x_2 - x_1) \\ \dot{x}_2 = (c + w_3)(x_1 + x_2) - (a - w_1)x_1 - x_1 x_3\,, \\ \dot{x}_3 = x_1 x_2 - (b - w_2)x_3 \end{cases} \tag{29}$$

with the parameters $(a,\ b,\ c) = (35,\ 3,\ 28)$ and the disturbances $\omega := (\omega_1,\ \omega_2,\ \omega_3)^T$. Then

$$A_0 = \begin{pmatrix} -a & a & 0 \\ c - a & c & 0 \\ 0 & 0 & -b \end{pmatrix},\ f_0(x) = \begin{pmatrix} 0 \\ -x_1 x_3 \\ x_1 x_2 \end{pmatrix},\ D_0(x) = \begin{pmatrix} x_1 - x_2 & 0 & 0 \\ x_1 & 0 & x_1 + x_2 \\ 0 & x_3 & 0 \end{pmatrix}.$$

When $\omega \equiv 0$, the attractor and state evolution of (29) are shown in Fig. 2. *(In this paper, the numerical results are derived with Euler-Maruyama method that solving stochastic differential equations [20]. And in all simulation, the stepsize is set as* $dt = 0.0001\ \mathrm{s}$*, and* $(a,\ b,\ c) = (35,\ 3,\ 28)$*, and the initial states* $x_0 = (2\ 8\ 10)^T$*.)*

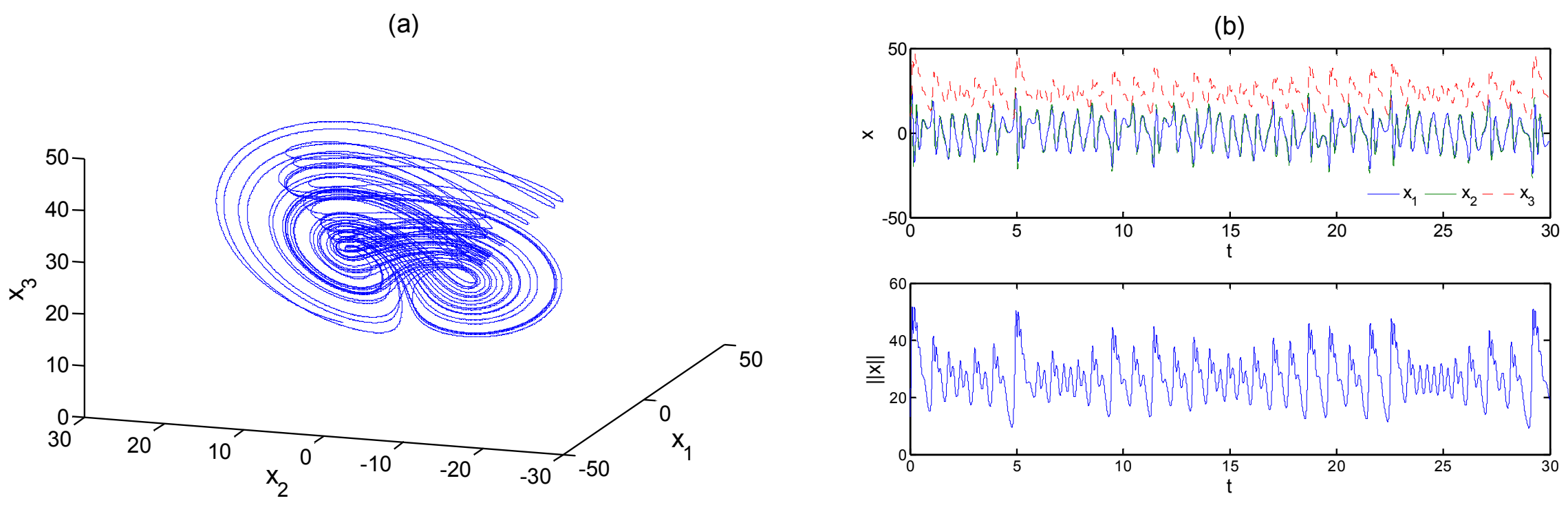

Fig. 2 (a) Chen attractor; (b) The corresponding state evolution.

Now, suppose the control system is

$$\begin{cases} \dot{x}_1 = (a + w_1)(x_2 - x_1) + [k_1(x) - ak_3(x)]x_1 + \sigma_1^c x_1 \dot{B}_1^c \\ \dot{x}_2 = (c + w_3)(x_1 + x_2) - (a + w_1)x_1 - x_1 x_3 + k_1(x)x_2 + ck_2(x) + \sigma_2^c x_2 \dot{B}_2^c , \\ \dot{x}_3 = x_1 x_2 - (b + w_2)x_3 + \left[\frac{2}{3}k_1(x) - bk_3(x)\right]x_3 + \sigma_3^c x_3 \dot{B}_3^c \end{cases} \quad (30)$$

where $k_i(x)$, $\sigma_i^c \dot{B}_i^c$ $(i = 1,\ 2,\ 3)$ are respectively the state feedback control law and the aiding noises. Then,

$$G(x) = \begin{pmatrix} x_1 & 0 & -ax_1 \\ x_2 & c & 0 \\ 2x_3/3 & 0 & -bx_3 \end{pmatrix},\ C(x) = diag(x_1,\ x_2,\ x_3),\ \sigma_c = diag(\sigma_1^c,\ \sigma_2^c,\ \sigma_3^c).$$

First, with $\sigma_c \equiv 0$, we are to find some state feedback controls $k(x) := [k_1(x) \quad k_2(x) \quad k_3(x)]^T$ so that system (30) is under the subthreshold region of being ISS w.r.t. $\omega$. And this should come after a designing task of stabilizing (30).

Since for (30) the assumption A3 stands, according to Lemma 1 we can equivalently stabilize system (29) to be robustly stable with a stability margin $\rho$ (a $\mathcal{K}_\infty$-function). Here, we set $\rho(\|x\|) = 0.5\|x\|$. Then if the disturbance inputs satisfy $\|\omega\| \leq \rho(\|x\|)$, according to Lemma 2 there should exists some GAS-CLF for (30). And we take such a candidate as $V(x) = 0.5\|x\|^2$. Further, let

$$\begin{aligned}
\dot{V}(x) &= x^T \dot{x} = x^T (A_0 x + f_0(x) + G(x)k(x) + D_0(x)\omega) \\
&\le x^T [A_0 x + f_0(x) + G(x)k(x)] + \rho(\| x \|) \| x^T D_0(x) \| ) \\
&= -a x_1^2 - b x_3^2 + c x_2 (x_1 + x_2) + [k_1(x) - a k_3(x)] x_1^2 + k_1(x) x_2^2 \\
&\quad + c k_2(x) x_2 + \left[ \frac{2}{3} k_1(x) - b k_3(x) \right] x_3^2 + \rho(\| x \|) \left\| \begin{bmatrix} -x_1^2 & -x_3^2 & x_2(x_1 + x_2) \end{bmatrix} \right\| \\
&\le -a x_1^2 - b x_3^2 + c x_2 (x_1 + x_2) + [k_1(x) - a k_3(x)] x_1^2 + k_1(x) x_2^2 \\
&\quad + c k_2(x) x_2 + \left[ \frac{2}{3} k_1(x) - b k_3(x) \right] x_3^2 + \rho(\| x \|) \left( \frac{3}{2} x_1^2 + \frac{3}{2} x_2^2 + x_3^2 \right) \\
&= 0,
\end{aligned}$$

with

$$k(x) = \begin{pmatrix} -1.5\rho(\| x \|) & -(x_1 + x_2) & -1 \end{pmatrix}^T . \tag{31}$$

Because $\dot{V}(x) = 0$ iff $x = 0$, $\dot{V}(x) < 0$ for all $x \neq 0$. Hence, $V(x)$ is a GAS-CLF for system (30). And thus, (30) is robustly stable. So system (30) is ISS w.r.t. $w$. This is confirmed by the trajectories shown in Fig. 3, with $\sigma_c \equiv 0$ and

$$\begin{aligned}
\omega_1 &= \sin(t) + 0.5\dot{B}_1^d, \\
\omega_2 &= 2\sin(t) + 0.25\dot{B}_2^d, \\
\omega_3 &= 0.5\sin(t) + \dot{B}_3^d,
\end{aligned}$$

where $\dot{B}_{1,2,3}^d$ denote the independent Gaussian white noises with distribution of $N(0,1)$.

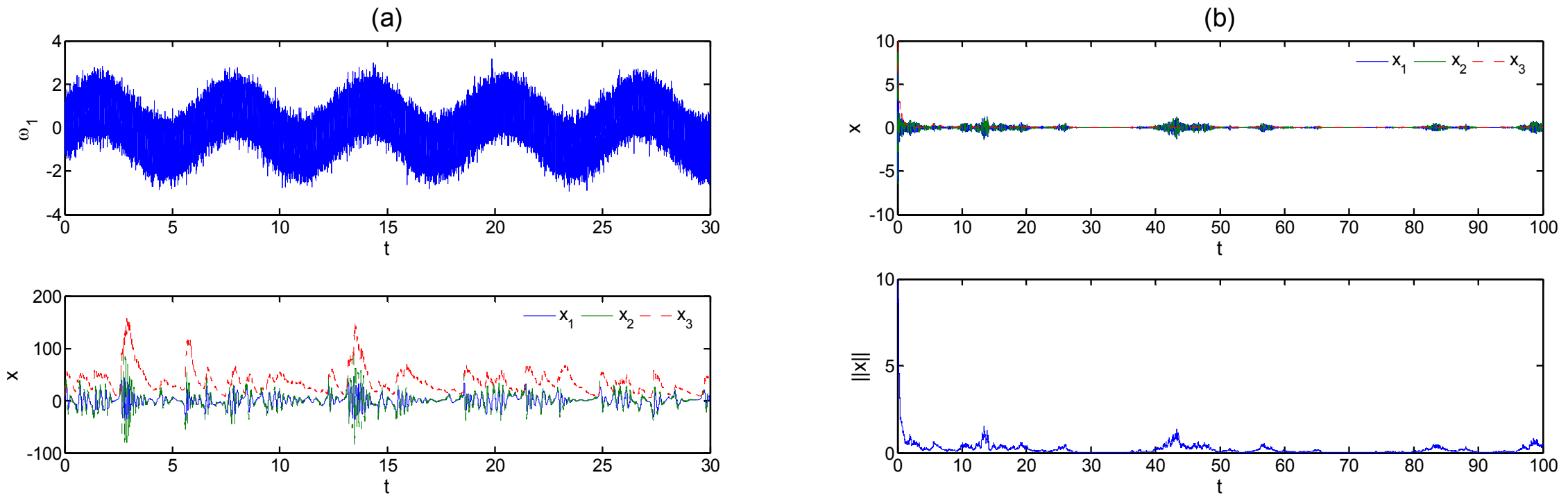

Fig. 3 (a) Noise $\omega_1$ and state evolution of system (30) with $k(x) \equiv 0$ and $\sigma_c \equiv 0$; (b) State evolution of system (30) with $\sigma_c \equiv 0$ and $k(x)$ taking values like (31).

However, if suppose $k(x)$ weakens to be

$$k(x) = \begin{pmatrix} -1.4\rho(\| x \|) & -0.9(x_1 + x_2) & -1 \end{pmatrix}^T , \tag{32}$$

then control system (30) specifically becomes

$$\begin{cases} \dot{x}_1 = ax_2 - 1.4\rho(\|x\|)x_1 + (x_2 - x_1)\omega_1 + \sigma_1^c x_1 \dot{B}_1^c \\ \dot{x}_2 = (0.1c - a)x_1 + 0.1cx_2 - x_1x_3 - 1.4\rho(\|x\|)x_2 \\ \qquad - x_1\omega_1 + (x_1 + x_2)\omega_3 + \sigma_2^c x_2 \dot{B}_2^c \\ \dot{x}_3 = x_1x_2 + 0.933\rho(\|x\|)x_3 - \omega_2 x_3 + \sigma_3^c x_3 \dot{B}_3^c \end{cases}. \tag{33}$$

And it will lose control, i.e., its solution will no longer keep close to the origin (see Fig. 4). And for system (33),

$$A = \begin{pmatrix} 0 & a & 0 \\ 0.1c - a & 0.1c & 0 \\ 0 & 0 & 0 \end{pmatrix},\ f(x) = \begin{pmatrix} -1.4\rho(\|x\|)x_1 \\ -x_1x_3 - 1.4\rho(\|x\|)x_2 \\ x_1x_2 + 0.933\rho(\|x\|)x_3 \end{pmatrix},$$

$$D(x) := \begin{pmatrix} D_1(x) & D_2(x) & D_3(x) \end{pmatrix} = D_0(x).$$

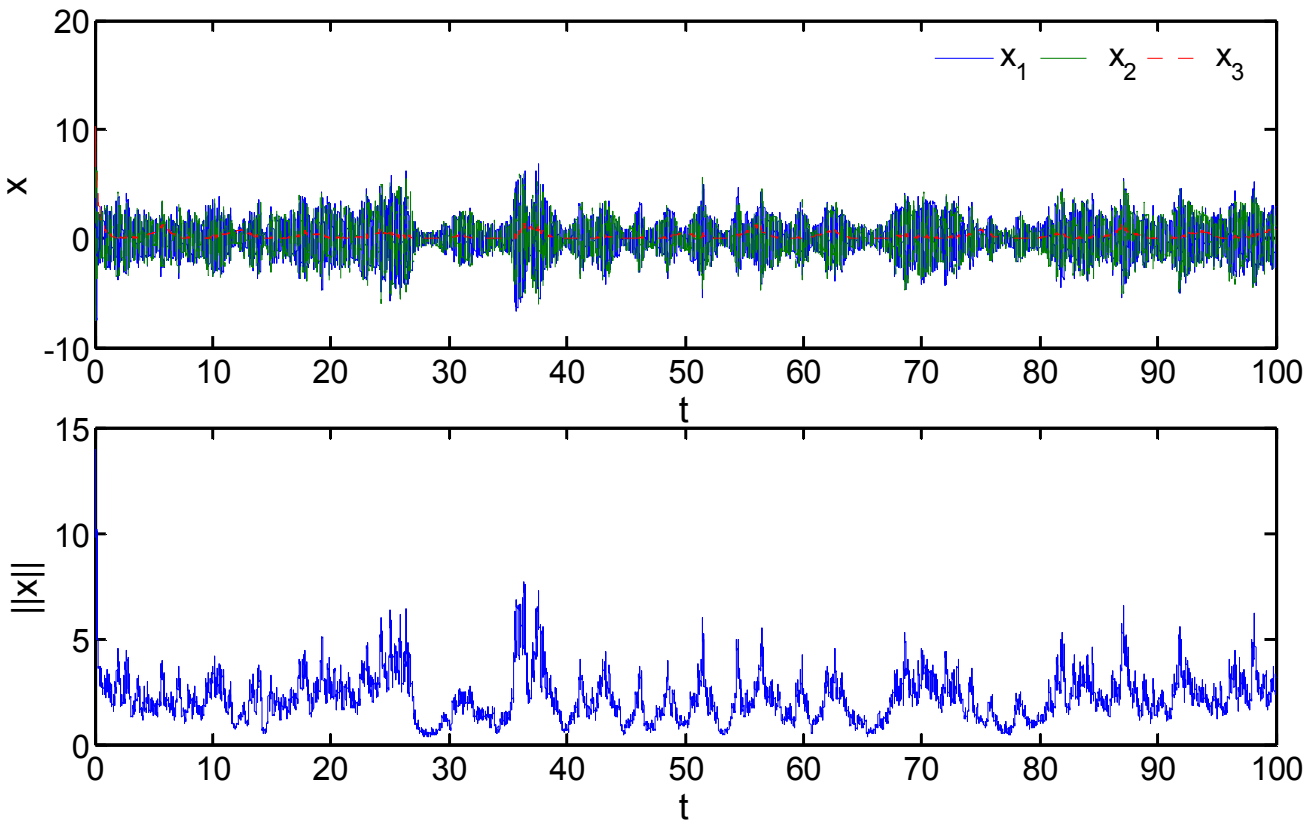

Fig. 4 Loss of control of system (30), with $k(x)$ being weakened to be (32).

From the knowledge of chaotic Chen system [15], we know that it is bounded when unperturbed. Here, we make an assumption that: system (33) keeps bounded under the disturbance $\omega$, namely it always holds that $\|x\| \le \theta$ for some nonnegative constant $\theta \in \mathbb{R}$. As a result, $f(x)$ should satisfy the Lipschitz condition, i.e., for some constant $L$ condition (10) holds.

Take the matrix $P$ in (15) as the unit matrix $I \in \mathbb{R}^{3\times 3}$. Then there exists such a matrix $R \ge 0$ that (12) hold, and it can be taken as $R = 0.75\theta I$. And then, according to Corollary 2, there would always exist such aiding noises that system (33) is ISS w.r.t. to $\omega$.

Set $\sigma_{1,2,3}^c = \sigma$. By numerical calculation, we obtain the deviation of the controlled dynamics from the origin (see Fig. 5. (a)), which is defined as (like that in [8])

$$\delta := \frac{\sum_{i=N_c}^{N} \| x_i \|^2}{N - N_c + 1},$$

where $N_c$ is the iterative time when the control is added and $N$ is the total iterative times. (In Fig. 5. (a), "Totally symmetric noise" means all the components of $\dot{B}_c$ and those of $\dot{B}_d$ are the same; "Symmetric noise" or "common noise" means all the components of $\dot{B}_c$ are the same and are independent of those of $\dot{B}_d$, whose components are independent of each other; "Independent noise" means all the components of $\dot{B}_c$ and those of $\dot{B}_d$ are independent of each other; "Asymmetric noise" means, in this example, $\dot{B}_2^c = \dot{B}_3^c = -\dot{B}_1^c$, and $\dot{B}_1^c$ and all the components of $\dot{B}_d$ are independent of each other.) And we find that when $\sigma \geq 3$ system (33) is ISS w.r.t. to $\omega$, in the common noise case. Take $\sigma = 3$ and the control effects are shown in Fig. 5. (b).

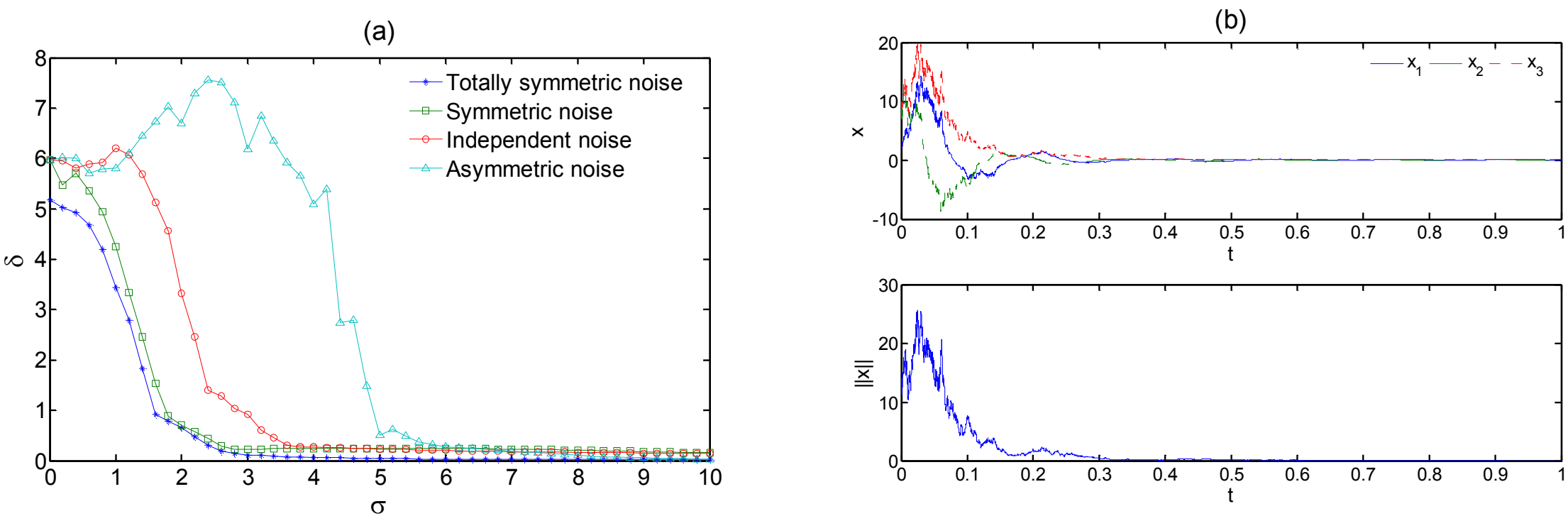

Fig. 5 (a) The deviation of the controlled dynamics from the origin as a function of the intensity of aiding noises, when $k(x)$ weakens to (32); (b) State evolution of system (33) with common noise as the aiding noise and $\sigma = 3$.

Moreover, the numerical results shows that even if the control $k(x)$ weakens to be

$$k(x) = \begin{pmatrix} -\rho(\| x \|) & -0.5(x_1 + x_2) & -1 \end{pmatrix}^T, \tag{34}$$

system (30) can be stabilized under the help of some proper white noises, see Fig. 6.

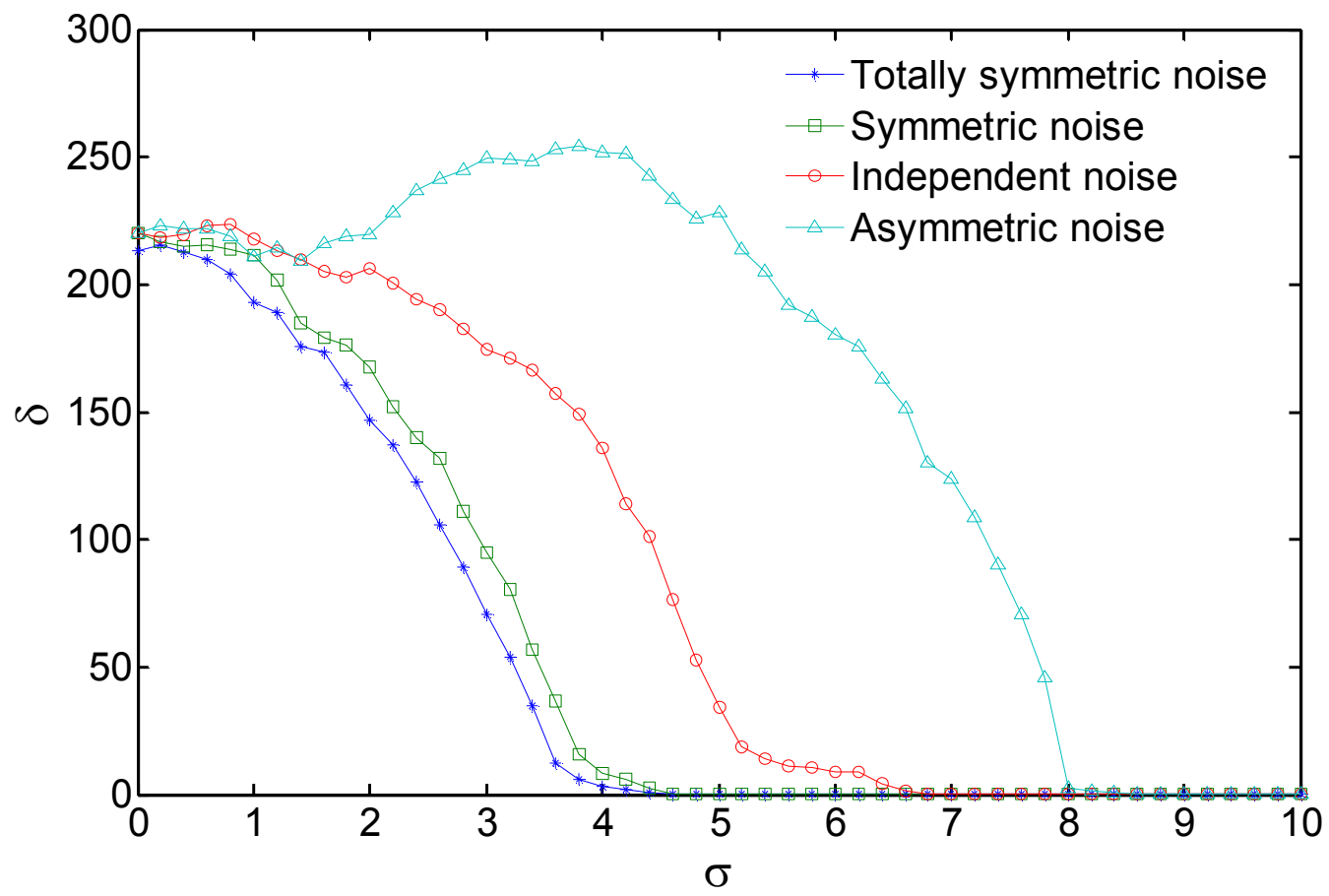

Fig. 6 The deviation of the controlled dynamics from the origin as a function of the intensity of aiding noise, when $k(x)$ weakens to (34).

So, in this sense, white noise can help to enhance the stability of systems even when the designed controls are under perturbation or not reliable enough.

Now, to show that aiding noise can also save the energy cost of the control, we define the cost function as

$$\psi = \frac{1}{t - t_c} \int_{t_c}^{t} \left( \left\| G(x(s)) k(x(s)) \right\|^2 + \left\| C(x(s)) \sigma_c \dot{B}_c(s) \right\|^2 \right) ds,$$

where $t_c$ is the time when the control is added and $t$ is the ending time of such integration. Take $t_c = 0$ and $t = 100\ s$. We calculate the cost of the un-noise-aided (i.e., $\sigma = 0$) case with controls being designed like (31), and the cost of the noise-aided case with controls being designed like (32) and aiding common noise with intensity $\sigma = 3$. And the cost are $\psi_1 \approx 2.35 \times 10^6$ and $\psi_2 \approx 1.41 \times 10^6$ respectively. Thus, the control cost is saved by using aiding noise. What is more, we can see from Fig. 5. (b) and Fig. 3. (b) that: in the noise-aided case the transient process is shortened greatly, comparing to the un-noise-aided case.

In addition, we investigate the impact of coherence of noises on the effect of control. And the cases of common (symmetric) noise and non-common (asymmetric or independent) noise are considered. The numerical results (see Fig. 6) show that, to guarantee the same control, the necessary intensity of common noise is the smallest and the independent noise larger and the asymmetric noise the largest. That is to say, common noise is the most suitable to aid control, and independent noise follows and last comes the asymmetric noise. This coincides with the conclusion derived in [9].

## 5. Conclusions and Discussions

In this paper, we theoretically prove the availability of utilizing white-noise to aid control—that is to help stabilize the perturbed system to be ISS w.r.t. the perturbations. Its effectiveness is also demonstrated by the given example, which confirms the previously observation that proper white noise can serve to reduce the control cost and shorten the transient process when realizing almost the same task control. Also, we should notice the general sense of the results obtained in this paper in explaining many noise enhanced or induced effects referred in the introduction.

However, we should also point out that there are still many aspects waiting to study in the noise-aroused effects. For one example, the noises are supposed to be independent of each other in this paper, but this may not keep true in reality. Although in some papers (e.g., [8]) there have reported that the coherence between noises seems having little influence on their aiding or enhancing effects, there are some countering observation that symmetric, independent and asymmetric noises contribute quite differently to the degree of system synchronization [9], which is also confirmed by the numerical results obtained here. With this being considered, the impact of coherent noises on stability or control is still in need of further theoretical analysis. For another example, the reasonable (or optimum) way of combining the state feedback and white-noise controls are yet unknown and much need to study in the future, as it seems unreasonable either to design costly feedback controls or to add aiding noises that are too small or too large. Besides, the effects of non-white noises on stabilization are not analyzed yet, which may have similar positive effects as white noise according to the results of paper [21]. Additionally, the recent articles [22-25] concerning stochastic stability should be alerted of, which may provide more general way to present the constructive role of noise in stability or control. So, in these senses, further systematic studies of constructive effects of noise on systems are still much in need.

## Acknowledgements

The author thanks X. Li, W. Lin, C. Cai, N. Xiao and D. Angeli for their helpful discussions on numerically solving stochastic differential equations. And also, the author greatly thanks G. Chen for his persistent concern and encouragement.